\documentclass[10pt,a4paper,reqno,article]{amsart}
\usepackage{amssymb,stmaryrd}
\usepackage{amsthm,amstext}
\usepackage{amssymb,amsmath}
\usepackage{amsfonts}
\usepackage{MnSymbol,slashed}
\usepackage{graphicx}
\usepackage{float}
\usepackage{latexsym}
\allowdisplaybreaks[1]
\usepackage{bbm}

\usepackage{dcolumn}
\newcolumntype{2}{D{.}{}{2.0}}

\usepackage{stmaryrd}

\theoremstyle{plain}
\newtheorem{theorem}{Theorem}[section] 

\newtheorem{lemma}[theorem]{Lemma}
\newtheorem{corollary}[theorem]{Corollary}
\newtheorem{proposition}[theorem]{Proposition}
\newtheorem{definition}[theorem]{Definition}
\newtheorem{example}[theorem]{Example}

\def\k{{{\mathbb{K}}}}

\newcommand{\CC}{\hbox{{$\mathcal C$}}}

\newcommand{\C}{\mathbb{C}}

\newcommand{\Z}{\mathbb{Z}}

\renewcommand{\k}{\mathbbm{k}}
\newcommand{\dimF}{\dim_1}

\newcommand{\sign}{\mathrm{sign}}

\newcommand{\isom}{{\cong}}
\newcommand{\eps}{{\epsilon}}

\newcommand{\tens}{\mathop{{\otimes}}}

\newcommand{\Ad}{\mathrm{ Ad}}

\newcommand{\id}{\mathrm{id}}

\renewcommand{\>}{\rangle}

\newcommand{\1}{1_A}

\newcommand{\prodf}{\includegraphics[scale=0.65]{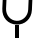}\ }
\newcommand{\coprodf}{\includegraphics[scale=0.7]{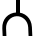}\ }
\newcommand{\und}[1]{\underline{#1}}
\newcommand{\Tr}{{\rm Tr}}
\newcommand{\loll}{\includegraphics[scale=0.65]{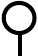}\ }
\newcommand{\uloll}{\includegraphics[scale=0.7]{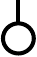}\ }

\begin{document}

\author{Shahn Majid${}^1$ and Konstanze Rietsch}
\address{School of Mathematical Sciences, Queen Mary University of London, Mile End Rd, E1 4NS \\ 
Department of Mathematics, Kings College London, The Strand}
\email{ s.majid@qmul.ac.uk, konstanze.rietsch@kcl.ac.uk}
\thanks{${}^1$Work done on sabbatical at Cambridge Quantum Computing}

\title{Planar spider theorem and asymmetric Frobenius algebras}
	\begin{abstract} The `spider theorem' for a general Frobenius algebra $A$, classifies all maps $A^{\otimes m}\to A^{\otimes n}$ that are built from the operations and, in a graphical representation, represented by a {\it connected} diagram. Here the algebra can be noncommutative and the Frobenius form can be asymmetric. We view this theorem as reducing any connected diagram to a standard form with $j$ beads $B$, where $j$ is the number of bounded connected components of the original diagram. We study the associated F-dimension Hilbert series $\dim_x=\sum_{j=0}^\infty x^j\dim_j$, where $\dim_j=\eps\circ B^j\circ 1$ are invariants of the Frobenius structure. We also study moduli of asymmetric quasispecial and `weakly symmetric' Frobenius structures and their  F-dimensions. Examples include general Frobenius structures on matrix algebras $A=M_d(\k)$ and on group algebras $\k G$ as well as on $u_q(sl_2)$ at low roots of unity.  \end{abstract}

\maketitle

\section{Introduction}

A Frobenius algebra $A$ is a unital algebra over a field $\k$ equipped with an invertible bilinear form $(\ ,\ )$ such that $(ab,c)=(a,bc)$ for all $a,b,c\in A$. By invertibility, we mean that the map $a\mapsto (a,\ )$ from $A\to A^*$ is invertible. This is equivalent to asking that $a\mapsto(\ ,a)$ is invertible, but note that these two maps may be different as we have not assumed that the bilinear form is symmetric.  In this setting the two inverse maps $A^*\to A$ can be encoded in a single element $g\in A\otimes A$ which is determined by the property that 
\[
\left((\ ,\ )\tens\id_A\right)(a\otimes g)=a=\left(\id_A\otimes(\ ,\ )\right)(g\otimes a)
\]
for all $a\in A$. The first equality says that $\varphi\mapsto (\varphi\tens\id_A)(g) $ inverts $a\mapsto (a,\ )$ while the second equality says that $\varphi\mapsto (\id_A\tens\varphi)(g)$ inverts $a\mapsto (\  ,a)$. Here, $\varphi$ denotes an element of $A^*$.

Let $\mu:A\otimes A\to A$ denote the product operation on $A$. We will call a Frobenius algebra {\em special} if $\mu$ sends $g$ to the identity element of $A$, 
\[ \mu(g)=\1.\]
Interpreting the bilinear form on $A$ as a linear map $(\ ,\ ): A\otimes A\to \k$ we may also evaluate $(\ ,\ )$ on $g$ to obtain a scalar,
\[ \dimF(A):=(\ ,\ )(g).\]
We call this the %`quantum' or 
\textit{Frobenius dimension} of $A$ or the F-dimension.

In this note, we study an invariant of Frobenius algebras motivated by a somewhat different, but equivalent, view of them  
as unital algebras $A$ which are also a  counital coalgebra in a compatible way\cite{Abr,CW}. Here, the additional structure is a  \textit{coproduct} $\Delta:A\to A\otimes A$ and a \textit{counit} $\eps:A\to \k$ obeying axioms dual to those of an algebra, and the further condition
\[ (\mu\tens\id_A)\circ(\id_A\tens\Delta)=\Delta\circ\mu=(\id_A\tens\mu)\circ(\Delta\tens\id_A).\]
From this is perspective $A$ is \textit{special}
if 
\[ \mu\circ \Delta=\id_A.\]
The equivalence between these two ways of thinking about Frobenius algebras is well-established: given a Frobenius algebra in the usual sense, as $(A,(\, ,\, ))$, use the bilinear form $(\, ,\, )$ to dualise the algebra to obtain a coalgebra. There is more than one way to describe the dualisation, but a left-right symmetric way is to characterise $\Delta(b):=\sum b_1\tens b_2$ (adopting a standard `Sweedler notation' for the coproduct) by
\[
\sum (a,b_1)(b_2,c)=(a,b c)=(ab,c)\]
for all $a, c\in A$. This is equivalent to setting $\Delta b=b g= gb$ in terms of the `inverse' $g$ of the Frobenius bilinear form. We also define $\eps(a)=(1_A,a)$ for all $a\in A$. Conversely, given the algebra-coalgebra version $(A,\mu,1_A,\Delta,\eps)$, define  $(\ ,\ )=\eps\circ\mu$ and $g=\Delta(1_A)$ to obtain a Frobenius algebra is the usual sense. From this point of view, we define the higher F-dimensions and F-Hilbert series
\[ \dim_j(A)=\eps(\mu(g)^j)=\eps\circ (\mu\circ \Delta)^j\circ 1,\quad \dim_x(A)=\sum_{j=0}^{\infty}x^j\dim_j(A)\]
where $\dim_1(A)=\dimF(A)$ is the previous quantum dimension. The equality of the first two expressions follows from the Frobenius properties. In the special Frobenius algebra case, clearly,
\[ \dim_x(A)={\dim_0(A)\over 1-x}.\]
Therefore the higher $F$-dimensions are of interest primarily in the non-special case. 

The definition of the higher $F$-dimensions is motivated from a so-called `spider theorem'. This theorem \cite{CD} says in the case of commutative $A$  that  any composition $A^{\tens m}\to A^{\tens n}$ of the above structure maps which corresponds to a connected graph in the diagrammatic notation is equal to a canonical map that looks like an $m-1$-fold product and an $n-1$-fold coproduct \begin{equation}\label{e:canonspider} 
(\id_{A}^{\tens n-2}\tens\Delta)\circ\cdots\circ (\id_A\tens\Delta)\circ\Delta\circ (\mu\circ\Delta)^j\circ \mu\circ(\id_A\tens\mu)\circ\cdots\circ(\id_{A}^{\tens m-2}\tens\mu).
\end{equation}
with some power $j$ of $\mu\circ\Delta$ in the middle. Here if $m=0$, the factors involving $\mu$ are replaced by $1$, while if $n=0$, the factors involving $\Delta$ are replaced by $\eps$. In these extreme cases the canonical form can be simplified further to  
$(\id_{A}^{\tens n-2}\tens\Delta)\circ\cdots\circ (\id_A\tens\Delta)\circ g$ or $(\ ,\ )\circ(\id_A\tens\mu)\circ\cdots \circ(\id_{A}^{\tens m-2}\tens\mu)$, respectively. A similar result holds for general Frobenius algebras provided we restrict to {\em planar} connected diagrams. We refer to this as the \textit{planar spider theorem} and a proof is provided in Section~\ref{secthm}, although we have since learned of at least one other in the recent literature\cite{Heu} (our proof is different and independent). In this case $k$ acquires the meaning of the number of bounded connected regions in the complement of the diagram, see Corollary~\ref{gen}. In the literature, there is a particular interest in the special Frobenius algebra case, and we prove the simpler, special version for this first. The non-special but symmetric case is also of interest and the spider theorem features in this case  implicitly and explicitly in \cite{LauPf,Q}. 

The planar spider theorem says in particular that any map $A^{\otimes m}\to A^{\otimes n}$ that is built from the structure maps as a connected planar diagram, in the case where $m=n=0$ must be multiplication by $\dim_j(A)$ for some $j$. For our theory to be useful we need a good supply of Frobenius algebras with asymmetric Frobenius forms. We provide explicit constructions of these on matrix algebras,  group algebras and the quantum group $u_q(sl_2)$ in Section~\ref{secex}, and study the F-dimensions. 

\subsection*{Acknowledgements}  SM would like to thank A. Kissinger for a helpful discussion motivating the project. 

\section{Diagram methods and planar spider theorem}\label{secthm}

We now introduce the diagrammatic notation as in \cite{Ma:alg,Ma:pri}. This notation was introduced in these works to do algebra in braided categories but is also useful in the monoidal case. Our results will apply to objects in any monoidal category, however the reader can keep in mind the category of vector spaces with unit object $\und 1=\k$, the ground field. In the diagrammatic notation,  morphisms are denoted by strands, possibly with labelled nodes. The identity morphism on a nontrivial object appears as a single vertical strand with no nodes. A tensor composite object is denoted by juxtaposition, and the product and coproduct in the case of an algebra and coalgebra appear respectively as  the trivalent vertices, \prodf and \coprodf. The unit object $\und 1$ is denoted by omission and the empty diagram is understood as the identity map $\und 1\to \und 1$. 

In our setting we will only be concerned with tensor composites of a single object $A$, which is to be understood at the external legs of all diagrams. Morphisms are read off the diagrams from top to bottom.  A diagram with $m$ strands coming in at the top and $n$ strands going out at the bottom represents a morphism from $A^{\otimes m}\to A^{\otimes n}$, keeping in mind that $A^{\otimes 0}=\und 1$.

The usual axioms of a special Frobenius algebra in the algebra-coalgebra form appear diagrammatically as 
\[ \includegraphics[scale=0.9]{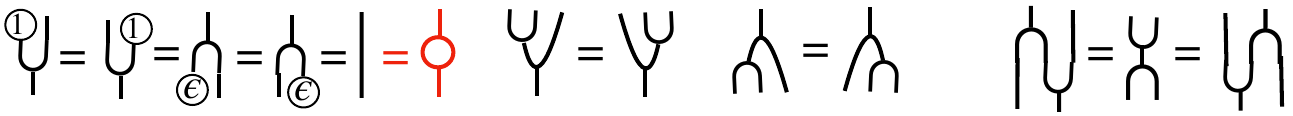} \]
where the unit $1=1_A$ is viewed as a map $1:\und 1\to A$ and the counit as a map $\eps:A\to \und 1$.  The last of the first group of equalities(in red) is the specialness condition. 

There are other possible characterisations of Frobenius algebras and one which is tailored to our purposes is given in the lemma below. Here, we group diagrammatic identities by the number $(m,n)$ of in and out strands.

\begin{lemma}\label{lem} Let $A$ be an object in a monoidal category. A Frobenius algebra structure on $A$ is equivalent to the specification of morphisms $\cup,\cap,\prodf,\coprodf$ and a morphism $1:\und 1\to A$ such that 
\begin{align}
\begin{aligned}[t]
\includegraphics{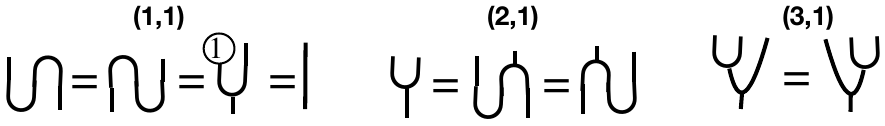}\label{d:coreaxioms}.
\end{aligned}
\end{align}
The further identities 
\begin{equation}\label{d:further}\includegraphics{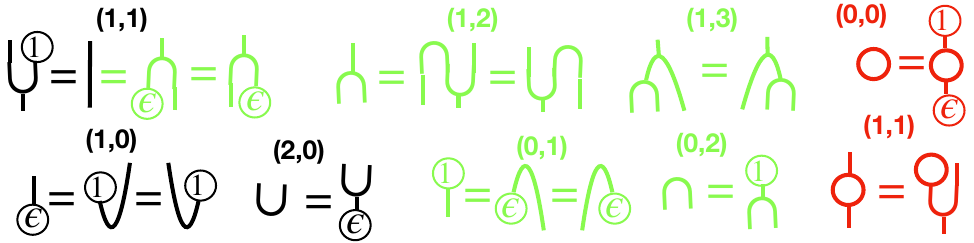}
\end{equation}
follow, with $\eps$ defined by the first  (1,0) identity. The Frobenius algebra $A$ in the form \eqref{d:coreaxioms} is \emph{special} if and only if $\loll=1$, or  equivalently if and only if $\uloll=\eps$.  
\end{lemma}

\proof To begin with, we assume the core identities \eqref{d:coreaxioms} and derive the others \eqref{d:further}. 

\begin{figure}
\[\includegraphics{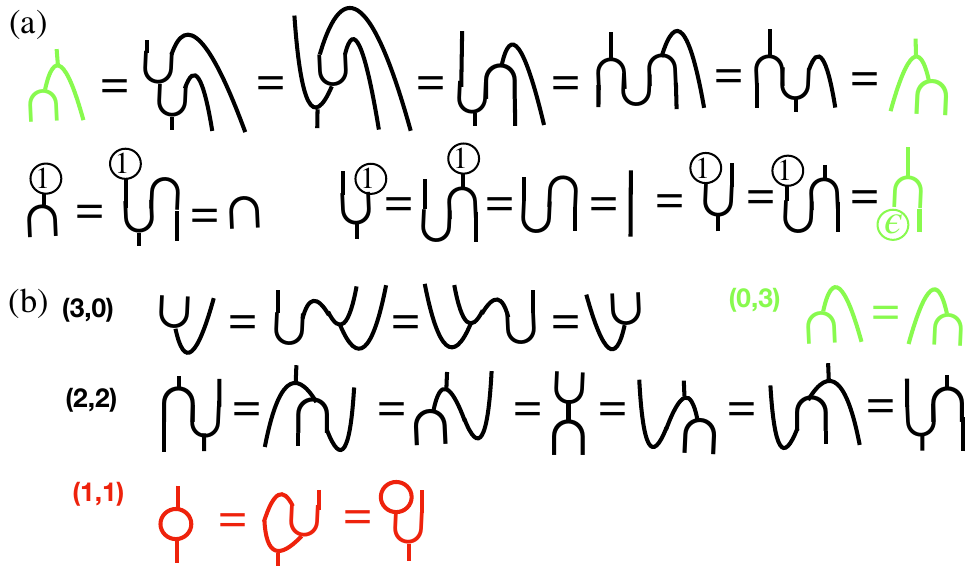}\]
\caption{\label{figlempf} Proof of Lemma~\ref{lem} with (a) derivation of redundant stated dual identities, (b) identities to obtain a Frobenius algebra and the special property.}
\end{figure}

\begin{enumerate}
\item[(1,2)] We connect the right output of $\cap$ to the left input of $\prodf$ and apply the first (2,1) identity to the $\prodf$ part of this diagram. Then by using the (1,1) assumption to replace the snake with the identity morphism,  the first equality in (1,2) follows. Similarly applying the left output of $\cap$ to the right input of $\prodf$ and then replacing $\prodf$ with the second form given in  (2,1) implies the second description of $\coprodf$ given  in (1,2). 
\item[(1,3)] Figure~\ref{figlempf} (a) uses the conversion between products and coproducts afforded by (1,2) and (2,1) to deduce the (1,3) identity from (3,1). 
\item[(0,2)] We replace the coproduct in $\coprodf\circ 1$ by the last expression in (1,2) and use the left identity assumption to deduce that $\Delta\circ 1=\cap$, as in Figure~\ref{figlempf}(a).
\item[(1,1)] Next, as shown in Figure~\ref{figlempf}(a) we compute $\prodf\circ(\id_A\tens 1)$ using the first identity of (2,1), and then apply (0,2) and the snake to prove the first of (1,1). 
\item[(1,1)] Next, as shown in the last part of Figure~\ref{figlempf}(a), we use $\prodf\circ(1\tens \id_A)=\id_A$ and the first identity of (2,1) again, to recognise the left counity property in (1,1) for $\eps$ defined as $\eps=\cup\circ (1\tens\id_A)$. 
\item[(2,0)] Now that we have proved that $\eps$ is a left counit, we can up-down reflect the proof of (0,2) to deduce the equality in (2,0).
\item[(1,1)] Now that we have shown (2,0), we up-down reflect  the proof in Figure~\ref{figlempf}(a) that $1$ is a right identity to prove that $\eps$ is a right counit, which is the last item in the main group of  (1,1) identities. 
\item[(1,0)] To deduce the second expression for $\eps$ in (1,0),  we use that $1$ is a right identity for the product. We compose this with $\eps$ in the only way possible and then simplify using (2,0) to obtain the result.  
\item[(0,1)] The (0,1) equalities follow by replacing $\eps$ using (1,0) and straightening out the diagram using (1,1) to get the unit.  
\item[(0,0)] the combination of the identities (0,2) and (2,0) implies the value  of the circle or F-dimension as stated in (0,0). 
\item [(1,1)] The last line in Figure~\ref{figlempf}(b) proves the extra (1,1) identity at the end of \eqref{d:further}.  We rewrite the coproduct using (1,2) and then use associativity (3,1). 
\end{enumerate}
We have now proved all of the further identities listed in the lemma. Let us finish the proof that $A$ is a Frobenius algebra, which we do by proving (3,0) and (2,2) in Figure~\ref{figlempf}(b).
\begin{enumerate} 
\item[(3,0)] As shown in Figure~\ref{figlempf}(b), we use the snake from (1,1) on either side and the second (1,2) identity in the middle. We have already established a unital algebra and $\cup$ is invertible via $\cap$, so we have a Frobenius algebra in the usual sense. 
\item [(2,2)] As shown in the second line of Figure~\ref{figlempf}(b), we rewrite the product using (2,1), use coassociativity (1,3) and then recognise $\Delta\circ\mu$ using (2,1). We then reverse these steps with left-right reflection to obtain the result. We have seen that $A$ is also a counital coalgebra, hence we have a Frobenius algebra in the algebra-coalgebra form. 
\end{enumerate}

Conversely, suppose we are given a Frobenius algebra in the algebra-coalgebra form. Then we let $\cup=\eps\circ\mu$ and $\cap=\Delta \circ 1$. Applying $\eps$ to the (2,2) axiom or evaluating it on $1$ immediately recovers the stated (2,1) and (1,2) identities respectively and doing the same to these gives the snake identities in (1,1).  

Equivalently, for completeness, suppose we are given a Frobenius algebra in the sense of a bilinear form $(\ ,\ )=\cup$ with inverse $g=\cap$, where the snake identities in (1,1) are assumed. We define the coproduct by the first of the  stated (1,2) identities. Equality of the outer parts of the first line of Figure~\ref{figlempf}(b) now tells us the other side of (1,2) and this converts by $\cup$  to (2,1) (an up-down reflection of the proof of (1,2) above). 

Hence the core identities in \eqref{d:further} are equivalent to a Frobenius algebra in the usual form or in the algebra-coalgebra form. The extra (1,1) identity in  \eqref{d:further} further implies that the property of being special is equivalent to saying that $\loll$ is a left identity, which is equivalent to $\loll=1$ given that the latter is already a two-sided identity. The same proof up-side down shows that the special property is also equivalent to $\uloll=\eps$. This picture of the \textit{special} property as $\loll=1$ is from \cite{Ma:zx}. \endproof

The last part of the lemma means that a \textit{special} Frobenius algebra can be thought of as an object $A$ together is a specification of morphisms $\cup,\cap,\prodf,\coprodf$ obeying the core identities, see \eqref{d:coreaxioms}, but replacing the unit map $1$ appearing in (1,1) by $\loll=\cap\circ\prodf$. The fact that $A$ is a \textit{unital} algebra follows by observing that $\loll$ has all of the properties of a unit element. 
 
 \begin{theorem}\label{thm}  A special Frobenius algebra is equivalent to the specification of morphisms $\cup,\cap,\prodf,\coprodf$ such that for each pair of nonnegative integers $(m,n)$, all compositions represented by a connected planar diagram with $m$ legs in and $n$ legs out are equal to the standard form \eqref{e:canonspider}, which diagrammatically appears as: 
  \[\includegraphics[scale=0.9]{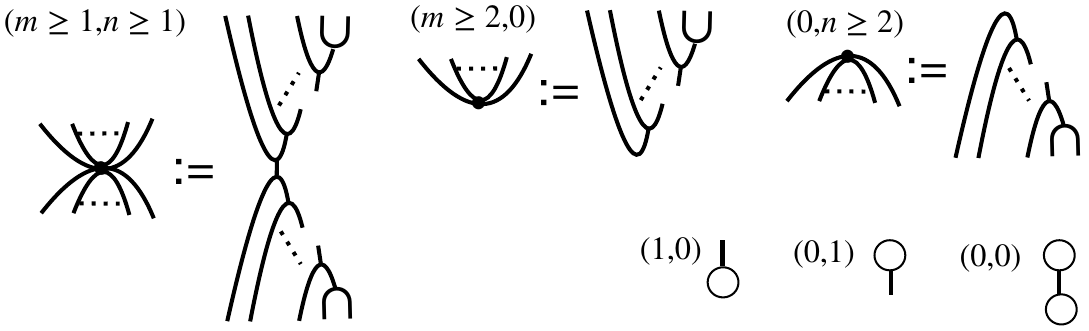} \]
 \end{theorem}
 \proof  
 To prove the theorem, we proceed by {\em induction on the total number of the operations} $\cup,\cap,\prodf,\coprodf$  in the connected planar diagram. For small numbers of operations, the statement of the theorem can be checked case by case using the identities in Lemma~\ref{lem}. For the induction step, given a connected planar diagram consisting of $n$ operations, we look at the last operation performed, adjusting the positioning if needed to obtain a unique one. Consider the diagram that remains after removing this last operation. It is either still a connected diagram or, if the last operation was either $\cup$ or $\prodf$, then it could be that it is the union of two connected diagrams. Either way the number of operations in each connected diagram is now less than $n$ and the induction hypothesis can be applied. 
 
\begin{figure}
\[\includegraphics{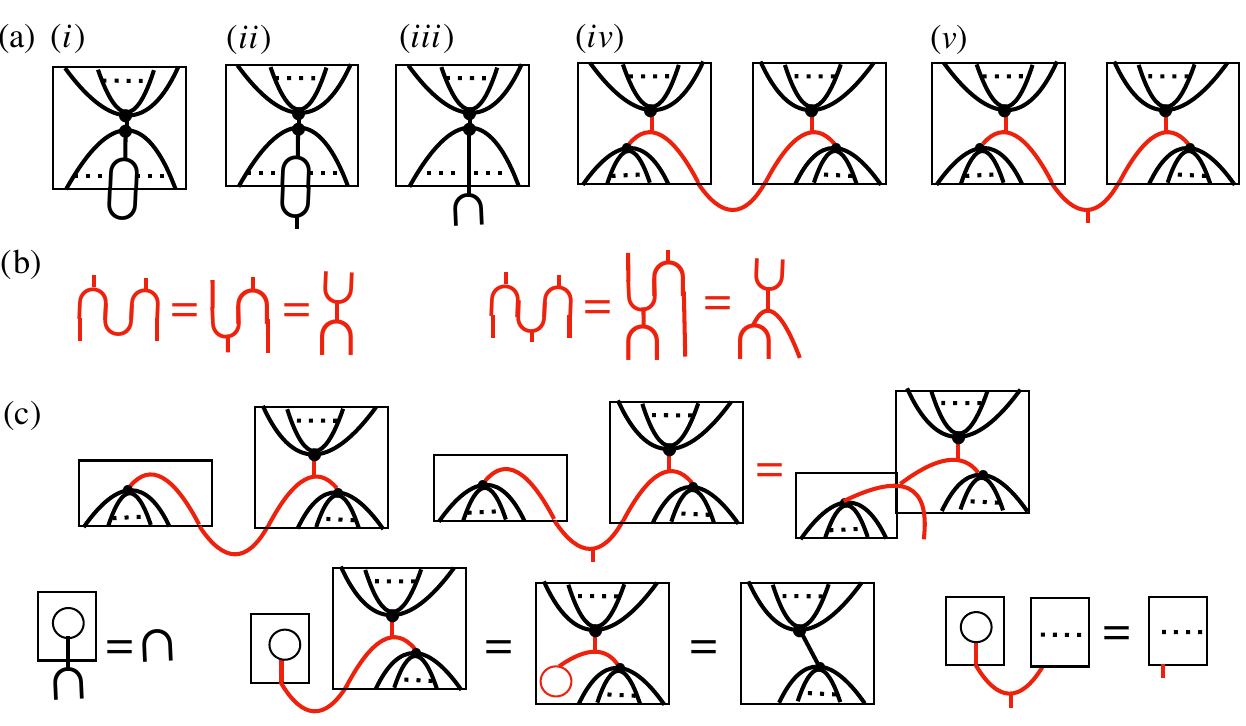}\]
\caption{\label{figthmpf} Proof of Theorem~\ref{thm} with (a) the combinations for the final operation in the generic case (b) identities in the proof and (c) special cases with $m=0$ on one side (or the other).}
\end{figure}

We thus obtain five possible scenarios which generically look like those shown in Figure~\ref{figthmpf} (a). The boxes contain connected diagrams with fewer than $n$ operations, and using the induction hypothesis we have already put them in the standard form.  In case (i) the whole diagram can be seen to be equivalent to one in the standard form by using the (2,0) identity from displayed equation \eqref{d:further}, with $\eps$ replaced by $\uloll$ as we are in the special Frobenius case. For (ii) the argument is similar, but using the usual definition of special. In case (iii) the diagram is seen to be equivalent to one in the standard form by using coassociativity. 
 
 We are now left with cases (iv) and (v). The red connecting diagrams in (iv) and (v) are reproduced in Figure~\ref{figthmpf} (b) and transformed using the (2,1) or (2,2) identities from Lemma~\ref{lem}. Replacing each red connecting diagram by its equivalent version from (b) gives a new diagram that can be transformed into the standard form using repeated applications of coassociativity. 
 
This covers most cases of the theorem. Let us check extreme cases, those that have either one or no output strand in one of the boxes. If there is only one output strand in either box, the proof just becomes even simpler using (3,0) or (2,1), associativity or the (2,2) identities from Lemma~\ref{lem}. This is left to the reader. We now turn to the case where one box has no input strands. If there are multiple strands coming out of this box and the other box is generic, then we are in the case shown in Figure~2 (c), top row. In this case the snake identity in (1,1) or a (1,2) identity can be used to put our diagram in standard form. If box has just one output strand so that is equal to $\loll$, then we must be in a version of case (iii), (iv) or (v). If we are in case (iii), we use the (0,2) identity in Lemma~\ref{lem}. In case (iv), we use the (1,0) identity in Lemma~\ref{lem} and the couninty property, where $\eps$ is replaced by $\uloll$ since we are in the special Frobenius algebra case.  In case (v) we use the unit property of $\loll$, for anything in the other box. The mirror images also hold. The case where there is no output strand does not occur, since the last operation is by definition outside of the boxes. \endproof

We have focused on the more elegant special Frobenius case, but we have structured our proof so that it applies also in the general non-special case with minor modifications. 

\begin{figure}
\[\includegraphics[scale=0.95]{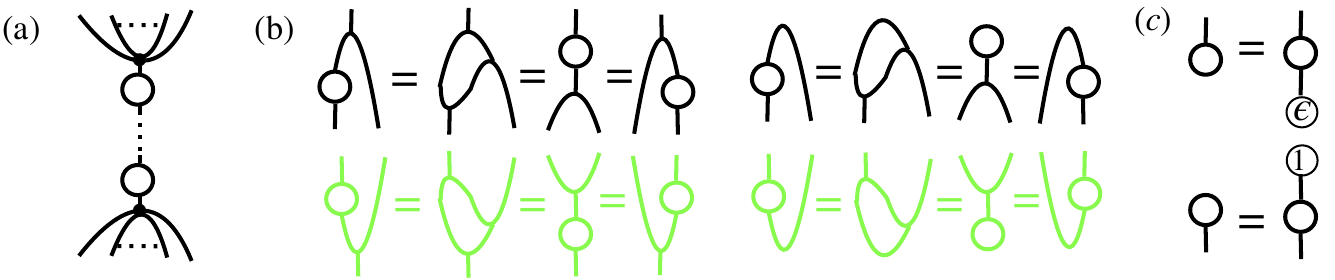}\]
\caption{\label{figgen} Proof of Corollary~\ref{gen} for general Frobenius algebras. Connected diagram compositions now have the standard form in (a) with `beads' studied in part (b). They can be created during the inductive step via (c).}
\end{figure}

\begin{definition}
 Let $j\ge 0$. For integers $m,n>0$ we define a `standard' diagram associated to $(m,n;j)$ as shown in Figure~\ref{figgen} (a), with $m$ strands going in, $n$ strands going out, and $j$ circle `beads' in the middle. If either $m$ or $n$ is equal to $0$ we define the standard diagram to have $1$ in place of the products, or $\eps$ in place of the coproducts, respectively. If both $m=n=0$ then we let the diagram given by $1$ followed by $j$ bubbles followed by $\eps$ be the standard diagram. 
\end{definition}

\begin{corollary}\label{gen}
In a Frobenius algebra $A$, any map $A^{\tens m}\to A^{\tens n}$ composed of  the operators $\cup,\cap,\prodf,\coprodf$ and $1,\eps$ that is represented by a \textit{connected} planar diagram is equal to the map given by the standard diagram with $m$ strands going in, $n$ strands going out, and $j$ beads, where $j$ is the number of bounded connected components of the complement of the original diagram. 
% with a fixed number $m$ of legs in and $n$ of legs out, 
\end{corollary}
\proof  The corollary will follow if we can show that any connected diagram composed of the operators $\cup,\cap,\prodf,\coprodf$ and $1,\eps$ is equivalent to one of the standard diagrams using the relations from Lemma~\ref{lem}. Notice that then the number of beads in the standard diagram is equal to the number of bounded connected components of the complement of the original diagram, since this number is an invariant for each of the transformations that we allow ourselves to apply to the diagrams.   

The first difference compared with proof of the theorem is that we could have morphisms $\eps,1$ in the mix. However, any $\eps$ in the composite is either at the start (in which case, for a connected graph, this is the whole graph and we are done) or connects to $1$ (in which case we are again done for a connected graph), or is applied to $\coprodf$ in which case we can remove it by the counity axiom, or is applied to $\prodf$ in which case we replace that by $\cup$ according to Lemma~\ref{lem}, or is applied to $\cap$ in which case  it gets replaced by $1$ according to Lemma~\ref{lem}. We then argue the same for $1$, which, unless it is the whole graph, must cancel with a product, convert a coproduct into a $\cap$ or connect to a $\cup$ and turn into an $\eps$. It is possible to have a chain of cups and caps through which this process repeats but this chain eventually terminates in a product or coproduct (so that the $1$ or $\eps$ disappear) or it brings a $1$ or $\eps$ to the start or end. Thus, it suffices to deal with  compositions of $\cup,\cap,\prodf,\coprodf$ other than the $(1,0),(0,1),(0,0)$ cases of $(m,n)$ which involve $\eps,1$ on one or both legs. 

We now follow the steps of the proof of the theorem. The key observation is in Figure~\ref{figgen}(b) that the `beads' can be taken through the products and coproducts at will. Hence they can always be collected in the middle if we are otherwise in the previous standard form. Referring to Figure~\ref{figthmpf}, the $\uloll$ created in (i) is now $\eps$ preceded by a `bead' as shown in Figure~\ref{figgen}(c). This, and likewise the `bead' in (ii) move up to join any bead already present in the box. The proof of the $m=0$ or $n=0$ cases in Figure~\ref{figthmpf}(c) is as before.   \endproof

In this context, we now define the higher F-dimensions of a Frobenius algebra as
\[ \includegraphics[scale=0.7]{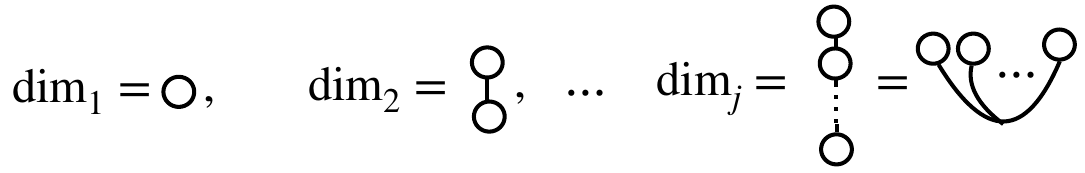}\]
where the diagram for $\dim_j$ has $j$ circles. We note that $\dim_j$ on the one hand is equal to $\eps\circ(\mu\circ\Delta)^j\circ 1$ and on the other hand, as shown in Figure~\ref{figlempf}(b) can also be written as $\eps\circ (\loll)^j$. Indeed this identity is an application of Corollary~\ref{gen}.

 \section{Examples of asymmetric Frobenius algebras and F-dimensions}\label{secex}

Asymmetric  noncommutative special Frobenius algebras do not seem to have been discussed particularly in the literature and it is an open question, at least in characteristic zero, as to whether the generality of our planar spider theorem is justified by significant classes of examples. We therefore conclude with a couple of constructions of them. We also consider more general examples of Frobenius algebras and their F-dimensions of interest in the asymmetric case. 

It can be convenient to work with the condition for a special Frobenius algebra only holding up to a nonzero scale factor, and we use the following definition. 

\begin{definition} A Frobenius algebra $A$ is called quasispecial if $\loll=\lambda 1$ for a nonzero scalar $\lambda$. In this case we call $\lambda$  the quasispecial scale factor. 
\end{definition}

Note that we require $\lambda$ to be nonzero to call the Frobenius algebra quasispecial, as one can always then change the normalisation of $(\ ,\ )$ (or the normaisation of the counit and, inversely, the coproduct) to make such a quasispecial Frobenius algebra special. 
%Namely this is just a matter of the normalisation of $(\ ,\ )$. 

Another natural scalar, this one associated to any Frobenius algebra, is 
\[\lambda':=\dim_0(A)=\eps\circ 1.
\]
This scalar is also in its own way dependent on the normalisation of $(\ ,\ )$ or equivalently of the Frobenius form $\eps=(1,\ )$, and we may refer to it as the \textit{counit scale factor}. For a quasispecial Frobenius algebra we note that the scale factors $\lambda$ and $\lambda'$ determine all of the F-dimensions. For example the primary F-dimension is 
\begin{equation}\label{e:lambdaprime}\dim_1(A)=\bigcircle=\eps(\loll)=\lambda\lambda',
\end{equation}
and this scalar is in fact canonical and independent of the normalisation. The linear form $\uloll$ is also independent of the normalisation for the same reasons, and indeed the F-dimension is equivalently expressed as $\uloll$ applied to $1$. We can immediately write down the F-Hilbert series in the quasispecial case as
\begin{equation}\label{quasidim} \dim_x(A)={\lambda'\over 1- \lambda x}.\end{equation}

A specific special or quasispecial, symmetric Frobenius structure, if it exists, will often be a convenient starting point for exploring all of the (symmetric and asymmetric) Frobenius structures on a given algebra $A$. To construct further examples, as well as to analyse moduli of other classes of Frobenius structures, we make use of a known observation that any two Frobenius structures on the same algebra are necessarily related by  the following `twisting' construction that we now recall. Let $(A,(\ ,\ ),g)$ be a Frobenius algebra and $u\in U(A)$ the group of invertible elements of $A$, and write $g=\sum g^1\tens g^2$ as an explicit notation. Then one can check that 
\[ (a,b)_u:= (ua,b),\quad g_u= \sum g^1\tens u^{-1}g^2\]
is another Frobenius algebra structure on $A$.  In the other direction, given two Frobenius structures $(\ ,\ )$ and $(\ ,\ )'$, use the former to identify $A\isom A^*$ as vector spaces by $a\mapsto (a,\ )$. The Frobenius linear form $\eps'=(1,\ )'\in A^*$  then corresponds to the required element $u\in A$ and has inverse given by $\eps$ applied to the first leg of the associated $g'$. 

The asymmetry of a Frobenius structure is measured by the Nakayama automorphism $\zeta$ of $A$, see \cite{Nak}, which is defined by $(\zeta(a),b)=(b,a)$. This automorphism changes to $\zeta_u=Ad_{u^{-1}}\circ \zeta$ in the twisted version, that is, 
\[ \zeta_u(a)=u^{-1}\zeta(a)u\]
is the Nakayama automorphism of the twisted Frobenius structure $(\ ,\ )_u$.
Similarly, $\loll$ and $\uloll$ are changed under twisting, to
\begin{equation}\label{e:lolltwist} \loll{\kern-3pt}^u= \sum g^1u^{-1}g^2,\quad\text{and}\quad \uloll_{\kern-3pt u}=\eps(ug^1u^{-1}g^2(\ ))=\eps(u\loll{\kern-3pt}^u(\ )).
\end{equation} 
This  means, in particular, that the property of being (quasi)special is not invariant  under twisting. Likewise, the F-dimension is not invariant under twisting. In fact, the F-dimension can be written as $\bigcircle=\Tr( \zeta)$ and therefore changes to $\bigcircle_u=\Tr(\Ad_{u^{-1}}\circ\zeta)$. Similarly for the higher F-dimensions. 

 We now introduce the following interesting subclass of Frobenius structures, which extends the class of symmetric ones.

\begin{definition} A Frobenius algebra is called {\em weakly symmetric} if $\uloll$ is a trace map, i.e. vanishes on $[A,A]$.
\end{definition}
Note that if a Frobenius algebra is symmetric then 
\begin{equation*} \uloll(ab)=(ab,\loll)=(a,b\loll)=(a,\loll b)=(\loll b,a)=(\loll, ba)=(ba,\loll)=\uloll(ba)\end{equation*}
for all $a,b\in A$, so symmetry implies weak symmetry. We used here that $\loll$ is central. The same argument proves the well-known result that a symmetric Frobenius form $\eps$ provides an identification $Z(A)\isom (A/[A,A])^*$ of the centre with the space of trace maps, via $c\mapsto \eps(c\,\underline{\ \ }\, )$. Thus, if we fix a symmetric Frobenius form as base for twisting,
then the twist by an invertible element $u$ is \textit{symmetric} if and only if $u\in Z(A)$, and it is 
\textit{weakly symmetric} if and only if $u\loll{\kern-3pt}^u\in Z(A)$, using \eqref{e:lolltwist}. %The latter condition is automatically satisfied if $u\in Z(A)$ since $\loll{\kern-3pt}^u$ is always central, but it is generally weaker since $\loll{\kern-3pt}^u$ need not be invertible. 

In other words, after fixing a base symmetric Frobenius form the space of all symmetric Frobenius forms can be identified with the group $Z(A)^{\times}$ of invertible elements of $Z(A)$; that is, it forms a $Z(A)^{\times}$-torsor. The action of $Z(A)^{\times}$ also extends to the space of weakly symmetric Frobenius forms,  which can be identified with the set of all invertible $u\in A$ such that $u\loll{\kern-3pt}^u\in Z(A)$. 

Weak symmetry in terms of $\loll$ is the assertion that $\loll$ is cocommutative, which amounts to the identity 
\begin{equation}\label{lollw} \loll(g-g_{21})=0\end{equation}
in $A^{\otimes 2}$, since $\loll$ is central. Here $g_{21}=\sum g^2\tens g^1$ in our notation. This again makes it clear that weak symmetry is implied by symmetry. 

From either point of view, it is clear that the class of asymmetric weakly symmetric Frobenius algebras is \textit{disjoint} from the class of  quasispecial Frobenius structures.

\subsection{Asymmetric forms on matrix algebras}
We have a standard symmetric quasispecial Frobenius algebra structure on the matrix algebra $A=M_d(\k)$ which is given by the pairing $\Tr(ab)$ for $a,b\in A$. Asymmetric Frobenius structures can be obtained from this one by twisting. 

\begin{proposition}\label{Mn} If  $u\in GL_d(\k)$ then $(a,b)=\Tr(uab)$ makes $A=M_d(\k)$ into a  Frobenius algebra which is {\em either} quasispecial with scale factor $\lambda=\Tr(u^{-1})$, when $\Tr(u^{-1})\ne 0$, {\em or otherwise}  asymmetric weakly symmetric, when $\Tr(u^{-1})=0$. The Frobenius structure is symmetric if and only if $u$ is a multiple of the identity. The lowest F-dimensions and Hilbert series are given by
 \[\lambda'=\dim_0(A)=\Tr(u),\quad {\textstyle\bigcirc}=\Tr(u)\Tr(u^{-1}),\quad 
 \dim_x(M_d(\k))= {\Tr(u)\over 1-\Tr(u^{-1})x}.\]
 Every Frobenius structure on $M_d(\k)$ is of this form. 
\end{proposition}
\proof That we have a Frobenius structure and that every Frobenius structure is of this form follow from the general remarks about twisting, but are worth seeing explictly. Here, $(ab,c)=(a,bc)$ is immediate and the bilinear form is clearly non-degenerate as $\Tr$ is nondegenerate as a bilinear form on $M_d(\k)$, and $u$ is invertible. We set $\cap=g=\sum_{i,j,k} E_{ij} \tens u^{-T}{}_{jk} E_{ki}$ where $u^{-T}$ denotes the transpose of $u^{-1}$ and $E_{ij}$ denotes the matrix with $1$ at row $i$ and column $j$ and zero elsewhere. Adopting a convention to sum over repeated indices and writing $a=a_{mp} E_{mp}$, we have
\[ (a,E_{ij})u^{-T}{}_{jk}E_{ki}=u_{lm}a_{mp} \delta_{pi}\delta_{jl}u^{-T}{}_{jk}E_{ki}=u_{jm}a_{mi}u^{-T}{}_{jk}E_{ki}=a_{ki}E_{ki}=a\]
\[ E_{ij}u^{-T}{}_{jk} (E_{ki},a)=E_{ij}u^{-T}{}_{jk} u_{lm}\delta_{mk}\delta_{pi}a_{pl}=E_{ij}u^{-T}{}_{jk} u_{lk}a_{il}=E_{ij}a_{ij}=a,\]
as required by the snake identities. Conversely, given any Frobenius structure, we can clearly write its  linear form, $\epsilon$, as $\Tr(u(\  ))$ for some invertible matrix $u$. 

Next, for the `lollipop', we have 
\begin{equation}\label{e:lollypopmatrix}\loll=\mu(g)=\sum u^{-T}{}_{jk} E_{ij} E_{ki}=\sum u^{-T}{}_{jk}\delta_{jk}E_{ii}=\Tr(u^{-1})1_A,\end{equation}
so this is  quasispecial with scale factor $\lambda=\Tr(u^{-1})$ provided this scalar is not zero. Clearly, $\eps(1)=(1_A,1_A)=\Tr(u)$, while the F-dimension is
\[ {\textstyle\bigcircle}=(\ ,\ )(g)=\Tr(u\mu(g))=\Tr(u)\Tr(u^{-1}).\]
For the bilinear form to obey $(a,b)=(b,a)$ for all $a,b$ requires $u^{-1}bu=b$ for all $b$, which implies that $u$ is a multiple of the identity and $\bigcircle=\dim(A)=d^2$ as expected. 
%More generally, the Nakayama automorphism is clearly conjugation by $u$, as our Frobenius structure is the twist by $u$ of the standard symmetric Frobenius structure. Further twist by $v\in GL_d(\k)$ simply changes $u$ to $uv$, $\eps(1)$ changes to $\Tr(uv)$ and the quasispecial scale factor changes to $\Tr(v^{-1}u^{-1})$. We see that the F-dimension is not invariant. 
Finally, we have that $\uloll=\Tr(u\loll(\, ))$, and the Frobenius structure is weakly symmetric if and only if $u\loll$ is central. Since $u\loll=\Tr(u^{-1})u$ this is the case precisely if  {\em either} $u$ is a nonzero multiple of the identity (which is the symmetric case) {\em or} if $\Tr(u^{-1})=0$. In the first case the multiple must be $d/\lambda$ so that $\bigcirc=d^2$, and the second case can only arise if $d>1$.

% given by evaluation against $\Tr$ of $u\Tr (u^{-1})$ which is central if {\em either} $u$ is a nonzero multiple  of the identity (which is the symmetric case and the multiple must be $d/\lambda$ so that $\bigcirc=d^2$) {\em or} $\Tr(u^{-1})=0$ which is possible if $d>1$. 
 \endproof

We see that the asymmetric weakly symmetric case gives rise to F-dimension~$0$. The converse is also true in the $d=2$ case, since then $\Tr(u^{-1})=\Tr(u)/\det(u)$. The following is not new but  entails an elementary constructive argument which we will then use further.  We need  $\k$ to be characteristic 0 and `sufficiently large' that the semisimple algebra of interest has a matrix block decomposition, for example algebraically closed.

\begin{lemma}\label{lemss} Let $A$ be a semisimple algebra over $\k$ of characteristic zero such that $A$ has a matrix  block decomposition. Then $A$ has a unique symmetric special Frobenius form.
\end{lemma} 
\proof We assume that $\k$ is such that $A$ is isomorphic to a direct sum of matrix algebras. Then by Proposition~\ref{Mn}, each matrix block has only one symmetric Frobenius form, up to scalar, namely standard one given by the trace (or a multiple of it). Each block by itself is a quasispecial Frobenius algebra, so the overall $\loll$ will be a sum of multiples of identity matrices. Therefore the Frobenius structure on $A$ is special if each of these summands is the identity in its respective block. This condition determines the Frobenius structure on $A$ uniquely, namely the Frobenius linear form  (which is the Frobenius counit) needs to have $\eps|_{M_d}=d\Tr$ on each $d\times d$ block. This is because scaling $\Tr$ by $d$ means scaling the metric $g$ and hence the $\loll$ for the block by $1/d$. Thus since the Frobenius linear form on $M_d(\k)$ given by $\Tr$ has quasispecial scale factor $d$, the form $\eps|_{M_d}=d\Tr$ will have $\loll=1_{M_d}$.    \endproof

The corresponding linear Frobenius form in this case plays a role similar to that of the integral on a Hopf algebra and provides a natural baseline for twisting. The general picture here, see \cite{Agu},  is that existence of a symmetric special Frobenius form corresponds to the algebra being `strongly separable' and in this case it is unique and necessarily provided by the trace form $(a,b)_{\Tr}=\Tr_AL_aL_b$ (where $L_a$ is left multiplication). Moreover, strong separability in characteristic zero is equivalent to semisimplicity. By contrast, a classic   result of Eilenberg and Nakayama is that an algebra being separable implies that it admits a symmetric Frobenius form. The converse is not true, unless the symmetric Frobenius form is special. By the same methods as in the proof of the above lemma, we see the following.

\begin{corollary} In the semisimple case of  $A=\oplus_{i=1}^{i=k} M_{d_i}(\k)$, a Frobenius structure obtained by twisting the unique special symmetric one by an invertible element $u=\oplus u_i$ has F-dimension and Hilbert series given by 
\[  \textstyle\bigcirc=\sum_{i=1}^{i=k}\Tr(u_i)\Tr(u_i^{-1}),\quad\text{and}\quad  \dim_x(A)=\sum_{i=1}^{i=k} {d_i\Tr(u_i)\over 1- d_i^{-1}\Tr(u_i^{-1})x},\]
respectively. The counit scale factor is $\lambda'=\sum_id_i\Tr(u_i)$. This Frobenius structure is quasispecial with scale factor $\lambda$ if and only if $\Tr(u_i^{-1})/d_i= \lambda$ for all $i$. 
\end{corollary}
\proof First, we compute $\dim_j(A)$. The special Frobenius structure has $\eps=\oplus_i d_i\Tr_{M_{d_i}}$, as in the proof of Lemma~\ref{lemss}. In the $i$-th block the twisting of $d_i\Tr_{M_{d_i}}$ by $u_i$ equals to the twist of $\Tr_{M_{d_i}}$ by $d_i u_i$. Thus \eqref{e:lollypopmatrix} from the proof of Proposition~\ref{Mn} implies that $\loll^u$ equals to $d_i^{-1}\Tr(u_i^{-1})I_{d_i}$ on the $i$-th block. Then
\[ \dim_j(A)=\eps_u((\loll^u)^j)=\eps_u\left(\oplus_i (d_i^{-1}\Tr(u_i^{-1})I_{d_i})^j\right)=\sum_{i=1}^k d_i\Tr(u_i)(d_i^{-1}\Tr(u_i^{-1}
))^j\]
using the twisted counit $\eps_u=\oplus_i d_i\Tr_{M_{d_i}}(u_i\,\underline{\ \ }\,)$. This then gives the stated formulas for $\dim_1$ and the F-Hilbert series. Setting $j=0$ gives the counit scale factor. Each block for which $\Tr(u_i^{-1})\ne 0$ is quasispecial, since $\loll^u=d_i^{-1}\Tr(u_i^{-1})I_{d_i}$ as we have seen. To be quasispecial overall, we just need these scale factors $d_i^{-1}\Tr(u_i^{-1})$ to agree for all $i$, in which case they give the quasispecial scale factor $\lambda$ as stated. We recover $\bigcirc=\lambda\lambda'$ and the F-Hilbert series (\ref{quasidim}) as before. \endproof 

It follows from this Corollary that there is more generally a $\dim(A)-k$ moduli space of special Frobenius structures parametrised by $k$-tuples of invertible matrices $u_i$ for which $\Tr(u_i^{-1})=d_i$ for each $i=1,\cdots,k$. The unique symmetric one among these has each $u_i$ equal to the identity. %a prescribed multiple of the identity. 

We now turn to the weakly symmetric Frobenius algebra structures on $A=\oplus_{i=1}^{i=k} M_{d_i}(\k)$. Recall that for any symmetric Frobenius algebra the F-dimension is always equal to $\dim(A)$, the usual dimension. 

\begin{corollary}\label{corwsss} In the semisimple case of  $A=\oplus_{i=1}^{i=k} M_{d_i}(\k)$, any weakly symmetric Frobenius structure has an integral F-dimension. Namely
\[ {\textstyle\bigcirc}=\sum_{i\in I}d_i^2\le \dim(A)\]
where $I\subseteq\{1,\cdots,k\}$ is some subset containing the set of all $i$ for which $d_i=1$. Such a weakly symmetric Frobenius structure is obtained by twisting the standard special Frobenius form by $u=\oplus u_i$ with  $\Tr(u_i^{-1})= 0$ if $i\notin I$, and $u_i=\mu_i^{-1}I_{d_i}$ otherwise, for some invertible $\mu_i\in \k$. The F-Hilbert series is then given by
\[ \dim_x(A)=\sum_{i\in I} {d_i^2\mu_i^{-1}\over  1- \mu_i x}.\]
\end{corollary}
\proof  The possible form of $u$ for weak symmetry to hold follows from Proposition~\ref{Mn} applied to each block. Namely, the Frobenius form is weakly symmetric on the $i$-th block if and only if $u_i\loll^{u_i}$ is in the centre of the block. This happens first of all whenever $\Tr(u_i^{-1})=0$ (so that $\loll^{u_i}=0$). Otherwise, if $\Tr(u_i^{-1})\ne 0$, Proposition~\ref{Mn} says that we are in the quasispecial case, and here the fact that $\loll^{u_i}$ is invertible combined with the condition to be weakly symmetric implies that $u_i$ itself must be in the centre of its block. Therefore $u_i=\mu_i^{-1}I_{d_i}$ for some nonzero element $\mu_i$ of the field. Then the quasispecial scale factor $\lambda_i$ on the $i$-th block is $\Tr(u_i^{-1})=\mu_i$ and the counit scale factor is $\lambda_i'=\eps_{u_i}(1)=d_i\Tr(\mu_i^{-1} I_{d_i} )=d_i^2\mu_i^{-1}$. Let $I$ be the set of all  $i\in\{1,\dotsc, k\}$ for which $\Tr(u_i^{-1})\ne 0$. We obtain a contribution of the form (\ref{quasidim}) to the F-Hilbert series for each such block, and correspondingly a summand of the form $\lambda_i\lambda_i'=d_i^2$  to the F-dimension. The other blocks, where $\Tr(u_i^{-1})=0$, do not contribute to the F-Hilbert series or the F-dimension. The stated formulas for F-dimension and F-Hilbert series follow. \endproof

It follows from this corollary that a weakly symmetric Frobenius form on $A=\oplus_{i=1}^{i=k} M_{d_i}(\k)$ is symmetric if and only if its F-dimension equals to $\dim(A)$, which is the $I=\{1,\cdots,k\}$ case above. And again we see that we have a $k$-dimensional moduli of symmetric Frobenius forms, with the unique special one at $\mu_1=\cdots=\mu_k=1$.  If a Frobenius form is both \textit{weakly symmetric} and \textit{special} then this also forces it to be the unique symmetric special one. 

\subsection{Asymmetric forms on group algebras}

Any finite dimensional Hopf algebra $H$ is a Frobenius algebra with Frobenius linear form $\eps$ given by the integral $\int$ of $H$, see~\cite{Par}. Here, $\int:H\to \k$ is defined by the property that $(\int\tens\id)\Delta_H=1\int$, where $\eps_H$ and $\Delta_H$ are the Hopf algebra counit and coproduct. (Note $\eps_H$ and $\Delta_H$ differ from $\eps$ and $\Delta$, the Frobenius counit and coproduct.) The inverse $g$ of the Frobenius algebra structure is then described in terms of the integral $\Lambda$ on its dual and the antipode $S$ of $H$, namely $g=(\id\otimes S)\Delta_H(\Lambda)$. 
% $\Lambda\in H$ is defined to obey $h\Lambda=\eps_H(h)\Lambda$ for all $h\in H$. 
This Frobenius structure is quasispecial with scale factor $\lambda=\eps_H(\Lambda)$, whenever this is nonzero. Otherwise, when $\eps_H(\Lambda)=0$, the F-dimension is zero. The counit scale factor $\lambda'=\int 1$ can also be zero, leading to F-dimension zero. Both scale factors are nonzero, or equivalently the F-dimension of $H$ is nonzero, if and only if $H$ is semisimple.

This construction of a Frobenius structure on $H$ can then be generalised by twisting as follows. Let  $u\in H$ be  any fixed, chosen invertible element, and use the notation $\Delta_H\Lambda=\Lambda_1\tens\Lambda_2$ (sum understood). %and $S$ denotes the Hopf algebra antipode.  
Setting $\eps(a):=\int ua$ determines a Frobenius structure with 
\[ (a,b)=\int u ab,\quad g=\Lambda_1\tens u^{-1}S\Lambda_2,\quad \loll=\Lambda_1 u^{-1}S\Lambda_2,\quad{\bigcircle}=\int u \Lambda_1 u^{-1}S\Lambda_2.\]
The counit scale factor is $\eps(1)=\int u$ and the associated linear map $\uloll$  on $H$ is
%where  $\eps$ denotes the Frobenius counit 
\[\uloll=\int u \Lambda_1 u^{-1}S\Lambda_2(\ ).\]

For a finite group algebra $H=\k G$, the integral is the characteristic function $\int=\delta_e$ extended linearly, where $e$ is the identity of the group $G$. The integral on the dual is $\Lambda=\sum_{g\in G} g$, and $\Delta_H(\Lambda)=\sum_{g\in G} g\otimes g$. Therefore  after twisting by  an invertible element $u\in \k G$ we have a Frobenius structure with
\[ \loll=\sum_{g\in G}gu^{-1}g^{-1}\quad\text{and}\quad{\textstyle\bigcircle}=\delta_e \big(\sum_{g\in G}ugu^{-1}g^{-1}\big).\] 
%That this is the most general form is also clear, by writing any Frobenius linear form as $\sum_g\mu_g\delta_{g^{-1}}=\delta_e(\sum_g \mu_g g(\ ))$ for some coefficients $\mu_g$. 
The Frobenius algebra  here is quasispecial if and only if %$\sum_{g\in \CC}gu^{-1}g^{-1}=0$ 
$\sum_{g\in G}gu^{-1}g^{-1}$ gives $0$ when projected onto any nontrivial conjugacy class $\CC$. In this case the quasispecial scale factor is given by $\lambda=|G|\delta_e(u^{-1})$, that is $\lambda=|G|\mu$ if $u^{-1}=\mu e+\cdots$ in the basis $G$. The counit scale factor $\lambda'$ is meanwhile determined by $u=\lambda'e+\cdots$ in the basis $G$. As a result the F-dimension is given by $\dim_1(\k G)=|G|\mu\lambda'=|G|\delta_e(u)\delta_e(u^{-1})$. The higher F-dimensions are 
\[ \dim_j(A)=\delta_e\big(\sum_{g_1,\cdots,g_j\in G} u g_1u^{-1}g_1^{-1}\cdots g_j u^{-1}g_j^{-1}\big).\]
Clearly, $\loll$ and all of the F-dimensions are invariant under conjugation of $u$ by elements of $G$, and hence so is $\dim_x(A)$. In particular, if we twist by $u\in G$ and associate to it the resulting F-Hilbert series, then this determines a power-series valued class function on $G$ which we also denote $\dim_x: G\to \k[[x]]$.

Note that the twisted Frobenius form is symmetric if $u\in Z(\k G)$, i.e. if it is an (invertible) linear combination of conjugacy classes viewed as elements of $\k G$. It is weakly symmetric, by definition, if 
\[ \uloll=\sum_{g\in G} \delta_e(ug u^{-1}g^{-1} (\  ))\]
vanishes on $[\k G,\k G]$, or equivalently, if it is a class function on $G$. That is, in the weakly symmetric case $\uloll$ is given by a linear combination of characters.

\begin{example}\label{exS3}\rm We let $G=S_3$ the symmetric group on 3 elements, and $r=(12),s=(23), t=rsr=(13)$ and $\k$ characteristic zero.
% and sufficiently large, for example algebraically closed. NOT NEEDED, all reps of S_n are defined over Z. 
We also recall the normalised characters $1,{\rm sgn},\chi$ of the trivial, sign and 2-dimensional representations, respectively (the latter has values $\{1,0,0,0,-{1\over 2},-{1\over 2}\}$ on the basis $\{e,r,s,t,rs,st\}$). We twist the standard Frobenius structure on $\k S_3$ by an invertible element $u$, though we may parametrise the twisted forms by $u^{-1}$. 

(a)  An element $u^{-1}$ of the required form for a \textit{special} Frobenius algebra is
\[ u^{-1}={e\over 6}+a(rs-sr)+b(r-t)+c(s-t)\]
with parameters $a,b,c\in \k$ satisfying
\[ d:=1+108(a^2-(b^2 + b c + c^2))\ne 0\]
for invertibility.
The inverse $u$, the twisting element, then looks like
\[ u={\textstyle\bigcircle} e+{6-{\textstyle\bigcircle}\over 2}(rs+sr)+9(2-{\textstyle\bigcircle})(a (rs-sr)+b(r-t)+c(s-t))\]
in terms of the F-dimension
\[{\textstyle\bigcircle}={2(2+d)\over d}.\]
In particular, we have a 3-parameter moduli of special Frobenius structures and these can have any F-dimension {\em other} than ${\textstyle\bigcircle}=2$. Of these there is a unique symmetric one, of dimension ${\textstyle\bigcircle}=6$, as only $u=e$ is in the centre. (This has to be the case by  Lemma~\ref{lemss}, given by the trace form.)  

(b) Looking beyond the `special' case, $Z(\k S_3)$ consists of elements of the form
\[ u=\alpha e + \beta(r+s+t)+\gamma(rs+sr)\]
for new parameters $\alpha,\beta,\gamma$. To obtain a symmetric Frobenius structure by twisting, we require this to be invertible, which it is if and only if
\[ (\alpha - \gamma)^4 (\alpha - 3 \beta + 2 \gamma) (\alpha + 3 \beta + 2 \gamma)\ne 0.\]
Hence there is more generally a 3-parameter moduli space of symmetric Frobenius forms, with just $u=e$ special among these.

(c) For the full moduli of all Frobenius forms, we take  general elements
\[ u=\alpha e + \beta(r+s+t)+\gamma(rs+sr)+ a (rs-sr)+b(r-t)+c(s-t)\]
which to be invertible needs
\[ (\alpha - 3 \beta + 2 \gamma) (\alpha + 3 \beta + 2 \gamma) ((\alpha - \gamma)^2 + 
   3 (a^2 -( b^2 + b c + c^2))^2\ne 0.\]
One can then compute $\uloll$ in these terms and ask for it to vanish on  $[\k S_3,\k S_3]$, i.e. on $r-t,s-t,rs-sr$, leading to two types of weakly symmetric Frobenius forms:  

(i) The 3-dimensional parameter space of symmetric forms in (b) where  
\[ a=b=c=0,\quad \uloll(\{e,r,s,t,rs,sr\})=\{6,0,0,0,0,0\};\quad \uloll=1+{\rm sgn}+4\chi.\]
(Note that $\uloll$ is the character of the regular representation as expected. 
 
 (ii) A 5-dimensional parameter space of non-symmetric weakly symmetric forms 
\[\alpha=\gamma,\quad \uloll(\{e,r,s,t,rs,sr\})=\{2,0,0,0,2,2\};\quad \uloll=1+{\rm sgn}.\]
This has no intersection with special forms (there are no asymmetric weakly symmetric special forms), which we can see here since the stated form of $u$ for the special case and $\alpha=\gamma$ would need ${\textstyle\bigcircle}=2$ there. 

(d) Finally, in the general parametrization in (c), the F-Hilbert series is
\[ \dim_x(A)=\frac{1}{6} \left(\frac{(\alpha+3 \beta+2 \gamma)^2}{\alpha+3 \beta+2 \gamma-6 x}+\frac{(\alpha-3 \beta+2 \gamma)^2}{\alpha-3 \beta+2 \gamma-6 x}+\frac{4 (\alpha-\gamma)^2}{\alpha-\gamma+{3\over 2}(2-\bigcircle) x}\right)\]
where 
\[ {\textstyle\bigcircle}=2+ {4(\alpha-\gamma)^2\over (\alpha - \gamma)^2 + 
   3 (a^2 -( b^2 + b c + c^2))}\]
is the F-dimension. This is a (non-linear) function on $\k S_3$ which is invariant under conjugation. Here $\alpha,\beta,\gamma,{\textstyle\bigcircle}$ are also such functions. When restricted to $u\in G$, we obtain a class function
\[ \dim_x={1\over 1-6 x}\left( {\rm A}+ {2x {\rm B}\over 1+ 6x}+ {3x {\rm C}\over 1+3x}\right).\]
where ${\rm A,B,C}$ are the characteristic functions of the trivial, order 3 and order 2 conjugacy classes respectively. This is not particularly convenient to write out in terms of the characters $1,\sign,\chi$. 

Also observe that $\bigcircle=6$ in the symmetric case where $a=b=c=0$ and for such Frobenius forms 
\[ \dim_x(A)=\frac{1}{6} \left(\frac{(\alpha+3 \beta+2 \gamma)^2}{\alpha+3 \beta+2 \gamma-6 x}+\frac{(\alpha-3 \beta+2 \gamma)^2}{\alpha-3 \beta+2 \gamma-6 x}+\frac{4 (\alpha-\gamma)^2}{\alpha-\gamma-6 x}\right).\]
By contrast, $\bigcirc=2$ for the asymmetric weakly symmetric case with $\alpha=\gamma$ and in this case
\[ \dim_x(A)=\frac{1}{2} \left(\frac{(\alpha+ \beta)^2}{\alpha+\beta-2 x}+\frac{(\alpha- \beta)^2}{\alpha-\beta-2 x}\right).\]
Note that the higher $F$-dimensions are not constant (independent of relevant parameters) even though $\bigcircle$ is in these two cases.  
\end{example}

These calculations were done directly using the group algebra. We can also use the block decomposition of $S_3$ as $\k\oplus \k\oplus M_2(\k)$ and more generally for any finite group $G$ to obtain such results.  
\begin{proposition} For any finite group $G$ and $\k$ characteristic zero and sufficiently large, for example algebraically closed, weakly symmetric Frobenius forms have 
\[ \uloll=\sum_{i\in I} d_i^2\chi_i,\quad \big|G/[G,G]\big|\le {\textstyle\bigcirc}=\sum_{i\in I}d_i^2\le |G|,\]
where  $d_i=\dim(V_i)$,  $\chi_i$ is the normalised character of the irreducible representation $V_i$ and  $I$ is a subset of irreducibles including all those of dimension 1.  For each $I$, the moduli space of weakly symmetric Frobenius forms is a product which has a $d_i^2-1$-dimensional factor for each $i\notin I$ and a $1$-dimensional factor for each $i\in I$.  
\end{proposition}
\proof The matrix blocks are assumed to be given by ${\rm End}(V_i)$ and on each block we have the description in Corollary~\ref{corwsss}. In the context of the Peter-Weyl decomposition of $\k G$, the trace in each block corresponds to the character $\chi_i$. Only those for which $\Tr(u_i^{-1})\ne 0$ are included in the sums, in which case $u_i$ is a multiple of the identity. Here $[G,G]$  denotes the commutator subgroup, $G/[G,G]$ is the Abelianisation of $G$ and $|G/[G,G]\big|$ coincides with the number of 1-dimensional irreducible representations. \endproof

In particular, the symmetric case corresponds to $I$ being the full set of all the irreducibles, with $\bigcirc=|G|$ maximal, and the moduli space of symmetric Frobenius forms has dimension the number of irreducible representations. For $\k S_n$ there are always $p(n)$ irreducible representations where $p(n)$ is the number of partitions of $n$, and only 2 of them (the trivial and the sign representation) have dimension 1. Hence, there are $2^{p(n)-2}$ strata to the moduli space of weakly symmetric Frobenius structures, with integer F-dimension
\[ 2\le {\textstyle\bigcirc}\le n!\]
and equality on the right exactly for symmetric Frobenius structures. Analysis of such moduli in the non-semisimple case is much harder and we end with a small example. 

\subsection{Asymmetric forms on $u_q(sl_2)$}

The reduced quantum groups $u_q(sl_2)$ over $\C$ are defined for all $q$ a primitive $n$th root of unity and have dimension $n^3$. We use the same conventions as in recent works \cite{AziMa,Ma:zx} with generators $K,E,F$ and relations
\[ K^n=1,\quad E^n=F^n=0,\quad EK=qKE, FK=q^{-1}KF,\quad [E,F]={K-K^{-1}\over q-q^{-1}}\]
where $[ , ]$ is the commutator. There is a quasitriangular Hopf algebra structure (triangular only for $n=2$) and it is known that the algebra is not semisimple. We take basis $\{K^iF^jE^k\}$. The integral is has $\int K F^{n-1}E^{n-1}=1$ and is zero on other basis elements. This together with the integral of the dual, $\Lambda=\sum_{i=0}^{n-1}K^iF^{n-1}E^{n-1}$ provides the standard Frobenius structure as discussed above for any finite-dimensional Hopf algebra. By the formulae there, we know that this has $\bigcirc=\loll=\uloll=0$ as the Hopf algebra is not semisimple. 

\subsubsection{The $n=2$ case.} We look first at the smallest in this family, namely $q=-1$.  In this case the generators $E,F,K$, have relations
\[ K^2=1,\quad E^2=F^2=[E,F]=\{E,K\}=\{F,K\}=0\]
where $\{\ ,\ \}$ is the anticommutator. We have $\int KEF=1$, and the integral is $0$ on the other $7$ elements of the monomial basis. The standard Frobenius structure has
\[ g=F\tens KE-KE\tens F+ E\tens KF-KF\tens E+ KEF\tens 1+1\tens KEF+ K\tens EF+ EF\tens K,\]
see \cite{Ma:zx}. We consider twists of this. 
 
A general invertible element here has the form 
\[ u=a 1+ b K + c EF + d KEF + \alpha E+ \beta KE + \gamma F + \delta KF;\quad  a^2\ne b^2\]
for coefficients in the field, with the one constraint shown. The inverse is then
\begin{align*} u^{-1}={1\over a^2-b^2}\Big(&a 1 - b K - { (a^2  + b^2)c - 2 a (b d + \alpha \gamma - \beta \delta) \over a^2-b^2}EF\\ &- {(a^2+b^2) d-2b(ac- \alpha\gamma+ \beta\delta)\over a^2-b^2}KEF - (\alpha E+ \beta KE + \gamma F + \delta KF)\Big).\end{align*}
From this data, we compute that the twisted Frobenius structure has linear form
\[ \eps_u=\{d,c,b,a,\delta,-\gamma,\beta,-\alpha\}\]
for the values on the basis in the same order as in $u$, and 
\[ \lambda'=\dim_0(A)=d,\quad \loll={8 b\over b^2-a^2}EF,\quad \uloll={8 b\over b^2-a^2}\{b,a,0,0,0,0,0,0\}\]
for the values again in the same basis order. Then
\[ {\textstyle\bigcirc}=\dim_1(A)={8 b^2\over b^2-a^2}\quad \text{ and }\quad \dim_x(A)=d + {8 b^2\over b^2-a^2}x,\]
as the higher $F$-dimensions all vanish. We see that there are {\em  no special Frobenius structures}. 

Next, from the product table, one can see that
\[ [u_{-1}(sl_2),u_{-1}(sl_2)]=\<KEF,E,KE,F,KF\>.\]
Hence we find that the entire 8-dimensional moduli of Frobenius structures are {\em all weakly symmetric} as $\uloll$ vanishes on commutators. Finally,  for a symmetric Frobenius structure, we need 
\[ a=\alpha=\beta=\gamma=\delta=0\]
for $\eps_u$ to vanish on commutators, hence there is a {\em 3-dimensional moduli} of these with parameters $b,c,d$ and $\bigcirc=8$ the usual dimension. These are further twists by an element of the centre (which is spanned by $1,EF,KEF$) of the natural symmetric choice coming from $u=K$.

\subsubsection{The $n=3$ case.} For $q=e^{2\pi\imath\over 3}$, we find that twisting the canonical Frobenius form obtained from the integral $\int$ by $u=K$ again gives a symmetric Frobenius form. The linear Frobenius form has $\eps_K(F^2E^2)=1$ and is zero on other basis monomials. We have the associated
\[ \loll=\frac{3 }{q}\left(c_q^2 - \frac{1}{ q^2}\right)\]
where 
\[ c_q=K+qK^{-1}-3q^2FE;\quad c_q^3=2+3q c_q\]
is the quadratic Casimir and one can check that it has the relation shown. Here $1,c_q,c_q^2,\Lambda=(1+K+K^2)F^2E^2$ are a basis of the centre. The  $F$-dimensions are
\[ \lambda'=\dim_0(A)=0,\quad {\textstyle\bigcirc}=27,\quad \dim_j(A)=3^{2j};\quad  j>1.\]
Any other symmetric Frobenius form can be obtained by a further twist of this by an invertible element of the centre, so $u=Kz$ where $z$ is in the centre. Hence we have a 4-dimensional moduli of symmetric Frobenius forms.

We now look for more general Frobenius structures. In the full 27 dimensional moduli of possible $u$, we focus first on those in the Cartan subalgebra, 
\[ u= u_0+u_1K+u_2K^2;\quad \delta=u_0^3+u_1^3+u_2^3-3 u_0u_1u_2\ne 0\]
with the condition shown for invertibility. Then
\[ \eps_u(F^2E^2)=u_1,\quad \eps_u(KF^2E^2)=u_0,\quad \eps_u(K^2F^2E^2)=u_2\]
for the linear form. The associated Frobenius form is not symmetric unless $u_0=u_2=0$. The $F$-dimensions are
\[ \lambda'=\dim_0(A)=0,\quad {\textstyle\bigcirc}=27 u_1{(u_1^2-u_0u_2)\over\delta},\quad \dim_j(A)=u_1\left(9{(u_1^2-u_0u_2)\over\delta}\right)^j.\]
In all cases $\loll$ is a polynomial of $c_q$ but never a nonzero multiple of $1$, and $\uloll$ does not vanish on
commutators except in the symmetric case, so there are no weakly symmetric forms in this class beyond the symmetric ones.

More generally, one can show that there are no special Frobenius structures at all on $u_q(sl_2)$ for $n=3$ (although, there is a 7-dimensional moduli where $\loll=0$). The key is to note that $\loll$ depends only on the terms of  $u^{-1}$  in the monomial basis with an equal power of $E,F$'s (as we will prove in general below), hence one only has to search in the 9-dimensional algebra of $u^{-1}$ with terms with equal powers of $E,F$.  Equally, we only need to search among $u$ of this type. One can also check that if $u=u_K+\cdots$ where  $u_K$ is a polynomial in $K$ and $\cdots$ are terms in the monomial basis with nonzero but equal powers of $E,F$'s and $u$ is invertible, then $u^{-1}=u_K^{-1}+\cdots$ and is of the same form. 

Also, there is a 4-dimensional moduli of such $u^{-1}$ for which the Frobenius form is not symmetric but the cosymmetry condition (\ref{lollw}) holds. Namely this consists of elements of the form
\[ u^{-1}= u_0(1+K)+  \frac{1}{3} (-q u_6 +(q+1) u_7- u_8)K^2 + \left(u_6+ K u_7+ K^2 u_8\right) F^2E^2  \]
for generic parameters $u_0,u_6,u_7,u_8$ (i.e., with some exclusions for invertibility). Hence, there is an at least 4-dimensional moduli of asymmetric weakly symmetric Frobenius forms. For example, 
\[ u=K(1-{3\over 2}F^2E^2), \quad u^{-1}=K^2+( - {3 q^2\over 2}+{(q-1)^2\over 2}K )F^2E^2 \]
leads to 
\[ \loll=6 (1 -  (q + 1)  K +  q K^2)FE + 18 F^2E^2,\]
\[ \uloll(1)=18,\quad \uloll(FE)=6,\quad \uloll(KFE)=-{18 q\over (2 + q)^2},\quad \uloll(K^2FE)={6  (q-1)\over 2 + q},\]
where $\uloll$ is $0$ on other basis elements. This Frobenius structure has 
\[ \dim_0(A)=-{3\over 2},\quad \bigcirc=18,\quad \dim_j(A)=0;\quad j>1,\]
and one can check that it is asymmetric weakly symmetric. 

\subsubsection{General $n$ case.} Motivated by our observations for $n=2,3$, we give some partial results in the general $n$ case.

\begin{proposition}\label{sl2sym} Twisting by $u=K$, i.e. 
\[ (h,g)=\int Khg\]
gives a symmetric Frobenius form.
\end{proposition}
\proof Using both the monomial basis and the reversed monomial basis, it suffices to show that
\[ (E^\alpha F^\beta K^\gamma, K^c F^b E^a)=\int q^{(\gamma+c)(\alpha-\beta)}K^{\gamma+c+1}E^\alpha F^{\beta+b}E^a=q^{(\gamma+c)(\alpha-\beta)}\delta_{\gamma+c,0}\delta_{\beta+b,n-1}\delta_{\alpha+a,n-1}   \]
and
\begin{multline*}
( K^c F^b E^a, E^\alpha F^\beta K^\gamma)= \int q^{\gamma(\alpha+a-\beta-b)}K^{c+\gamma+1}F^bE^{a+\alpha}F^{\beta}\\ 
=q^{\gamma(\alpha+a-\beta-b)}\delta_{c+\gamma,0}\delta_{b+\beta,n-1}\delta_{a+\alpha, n-1}  
\end{multline*}
coincide for all powers. In both cases we used the $K$ commutation relations to collect all powers of $K$ to the left and the $[E,F]$ relations repeatedly to normal order with the $F$'s to the left of the $E$'s. However, the normal ordered terms arising from $[E^\alpha,F^{\beta+b}]$ have lower powers of $F$ than $\beta+b\le n-1$ and likewise from $[E^{a+\alpha},F^\beta]$ have lower powers of $E$ than $a+\alpha\le n-1$, and hence cannot contribute to the integral which in the standard basis has support only on $K F^{n-1}E^{n-1}$. Moreover, the $q$ factor disappears in the  first case as $\gamma+c=0$ and in the second case as $a+\alpha=b+\beta$ (as both are equal $n-1$), due to the delta-functions. \endproof

This choice $u=K$ provides a natural base point; it follows that all other symmetric Frobenius forms are given by a further twisting by invertible elements of the centre, so the dimension of $Z(A)$ is the dimension of the moduli of symmetric Frobenius forms. For example, in  \cite{Ker} it is shown that the centre is $(3n-1)/2$-dimensional if $n$ is an odd prime.

Next, if we look at the class of forms given by  $u$ in the Cartan subalgebra,
\[ u= \sum_{i=0}^{n-1} u_i K^i,\quad u_i\in \C,\]
then invertibility of $u$ is equivalent to the circulant condition
\[ \det\begin{pmatrix}u_0 & u_{n-1}& u_{n-2}& \cdots & u_{2}& u_{1}\\ u_{1}& u_0& u_{n-1}&  \cdots &  u_{3} & u_{2}\\
\vdots & & &\cdots  & & \vdots\\ 
u_{n-2} & u_{n-3}   & u_{n-4} &    \cdots &u_0& u_{n-1}\\
u_{n-1} & u_{n-2}   &  u_{n-3}&  \cdots &u_{1}& u_0\end{pmatrix}\ne 0.\]
It can be expected on the basis of our $n=3$ result that there are no asymmetric weakly symmetric Frobenius structures in this class of $u$, as well as no special ones. 

More generally, we note that $u_q(sl_2)=u_q(sl_2)_0\oplus u_q(sl_2)_1\cdots \oplus u_q(sl_2)_{n-1}$ according to the the $\Z/n\Z$ grading induced by conjugation by $K$. Namely, the degree of a monomial is the number of $F$'s minus the number of $E$'s in our conventions. The degree $0$ part is the $n^2$-dimensional commutative subalgebra generated by $K$ and $FE$. Here $(FE)^n$ can be reduced to an element of the algebra with smaller powers of $FE$ in view of the commutation relations and $F^n=E^n=0$. Moreover, $Z(A)\subset u_q(sl_2)_0$ as a subalgebra since its elements must commute with $K$.

The following proposition implies that any element of $u_q(sl_2)$ that can arise as  $\loll$ from a choice of Frobenius structure, can in particular be obtained by twisting the standard Frobenius structure by some element of degree $0$.
\begin{proposition} Consider the standard Frobenius structure on $u_q(sl_2)$ given by $\int$ and its twisting by an invertible element $u$.  Then the twisted $\loll$ depends only on the part of $u^{-1}$ that lies in the commutative subalgebra $u_q(sl_2)_0$ under the grading decomposition.
\end{proposition}
\proof First note that $u_q(sl_2)_0$ is closed under inversion, so the statement is equivalent to saying that $\loll$ depends only on the component of $u^{-1}$ with equal numbers of $E,F$'s in the monomial basis. We start with an observation about the metric $g$ for the standard Frobenius structure. Let $\Lambda_K=\sum_{i=0}^{n-1} K^i$, so that $\Lambda=\Lambda_KF^{n-1}E^{n-1}$. Therefore we have 
$\Delta_H\Lambda=\Delta_H (\Lambda_K F^{n-1}E^{n-1})=(\Delta_H \Lambda_K)(\Delta_H F^{n-1})(\Delta_H E^{n-1})$. The Hopf algebra coproduct takes the form
\[ \Delta_H K=K\tens K,\quad \Delta_H E=1\tens E+E\tens K ,\quad \Delta_H F=F\tens 1+K^{-1}\tens F ,\]
on the generators. Therefore $\Delta_H(\Lambda_K)$, involves only $K$, and $\Delta_H(F^{n-1})$, has total power $F^{n-1}$ spread between the tensor factors, and no $E$'s, while $\Delta_H(E^{n-1})$ has a total power $E^{n-1}$ split between the tensor factors, with no $F$'s occurring. The antipode $S$ does not change the powers of $E$ or $F$ but does reverse the order. Hence the standard metric before twisting takes the form
\[ g= \sum\Lambda_1\otimes S \Lambda_2=\sum_{i,j,l,m}c_{ijlm} K^i F^l E^m \tens K^j E^{n-1-m}F^{n-1-l}\]
for some coefficients $c_{ijlm}$. The twisted metric is obtained by inserting $u^{-1}$ in front of the antipode, and the twisted $\loll$ is $\sum\Lambda_1 u^{-1} S \Lambda_2$. It follows that any terms in $u^{-1}$ with an unequal power of $E,F$ contribute terms in the twisted $\loll$ which also have an unequal power. But these terms must in the end all cancel since $\loll\in Z(A)$, which is contained in $u_q(sl_2)_0$, so that all of its terms have equal powers of $E,F$.  \endproof

Both the `special' property and the weakly symmetric property depend only on $\loll$, in the latter case via the cocommutativity property (\ref{lollw}). Hence, this proposition tells us that Frobenius forms of either type can be reached, if they exist, working only with $u^{-1}$, or equivalently $u$, in $u_q(sl_2)_0$.  Based on $n=3$, we expect significant moduli of asymmetric weakly symmetric Frobenius forms for all $n>2$, but no special Frobenius forms.  The actual analysis will not be attempted here, but we note that concerning inverses, our observation from $n=3$ holds in general:

\begin{lemma} If $u=u_K+ \cdots\in u_q(sl_2)_0$ is invertible then its inverse has the form $u^{-1}=u_K^{-1}+\cdots$ where $u_K$ is the Cartan subalgebra and $\cdots$ refers to terms with $E,F$ in the standard monomial basis. 
\end{lemma}
\proof Let $u=u_K+\sum_{\alpha=1}^{n-1}u_\alpha F^\alpha E^\alpha $ and $u^{-1}=v_K+\sum_{\beta=1}^{n-1}v_\beta F^\beta E^\beta$ with $u_K,v_K, u_\alpha,v_\beta$ in the Cartan subalgebra. The product is $uu^{-1}=u_K v_K+\sum_{\alpha,\beta}u_\alpha v_\beta F^\alpha E^\alpha F^\beta E^\beta$. Writing  $F^\alpha E^\alpha F^\beta E^\beta=F^{\alpha+\beta} E^{\alpha+\beta} + F^\alpha [E^\alpha, F^\beta] E^\beta$, these commutators when normal ordered appear with lower powers of $E,F$, as in the proof of Proposition~\ref{sl2sym}, but the contribution to $uu^{-1}$ still enters with positive powers of $E,F$. Hence we need $u_K v_K=1$. Similarly for $u^{-1}u$.\endproof

One can also look at these issues for the Taft algebra $u_q(b_+)\subseteq u_q(sl_2)$ generated by $K,F$ only in our conventions. This has $\int F^{n-1}=1$ as the only support in the monomial basis. At least for $n=3$, one can check that $\loll=\uloll=0$ for all Frobenius forms.

\end{document}